\documentclass{article}
\usepackage{amsmath,amsthm,amsfonts,amssymb}
\usepackage{mathrsfs}
\usepackage{soul}
\usepackage[colorlinks,linkcolor=blue,anchorcolor=blue,citecolor=blue]{hyperref}
\usepackage[numbers,sort&compress]{natbib}
\usepackage{changes}
\numberwithin{equation}{section}

\newtheorem{definition}{Definition}[section]
\newtheorem{theorem}[definition]{Theorem}
\newtheorem{proposition}[definition]{Proposition}

\newtheorem{lemma}[definition]{Lemma}
\newtheorem{corollary}[definition]{Corollary}

\title{On stability of exponentially subelliptic harmonic maps.}
\author{Xin Huang}
\date{}

\begin{document}
\maketitle
	
\begin{abstract}
In this paper, we study the stability problem of  exponentially subelliptic harmonic maps from sub-Riemannian manifolds to Riemannian manifolds. We derive the first and second variation formulas for exponentially subelliptic harmonic maps, and apply these formulas to prove that if the target manifold has nonpositive curvature, the exponentially subelliptic harmonic map is stable. Further, we obtain the instability of exponentially subelliptic harmonic maps when the target manifold is a sphere.

\hspace*{\fill} 
\end{abstract}

	\section{Introduction} 
	Exponentially harmonic maps were first introduced
	by Eells and Lemaire \cite{EL92} as a critical point of the exponential energy functional. Later, 
	Hong and Yang \cite{HY93}, Zhang, Wang and Liu \cite{ZWL98}, Cheung and Leung\cite{CL99}, Chiang et al. \cite{CH07, Ch15, Ch16} studied the stability of exponentially harmonic maps and related properties. Among them, Chiang and Yang \cite{CH07} proved that the exponentially harmonic map between two Riemannian is stable if the target manifold has non-positive sectional curvature. Furthermore, Chiang \cite{Ch16} proved that the non-constant exponentially harmonic
	map from a compact Riemannian manifold into
	a sphere $S^n$ $(n \ge 3)$ is unstable when $|df|^2 < n-2$. On the other hand, Barletta and Dragomir \cite{BD04} studied the pseudo harmonic maps (with potential) and their stability problem. Later, the author\cite{Don21} introduced general subelliptic harmonic maps from sub-Riemannian manifolds and the authors \cite{DLY22} generalized the results to subelliptic harmonic maps with potenital. Recently, Chong, Dong and Yang \cite{CDY22} considered the stability problem of subelliptic harmonic maps (with potential) from sub-Riemannian manifolds to another Riemannian manifolds and offered Leung type results when the target manifold is a sphere $S^n$ $(n \le 3)$.

	In this paper, we study the stability of exponentially subelliptic harmonic maps from
	a compact sub-Riemannian manifold into a Riemannian manifold. First, we derive the  first variation formulas (see Theorem \ref{th2.3}) for exponentially subelliptic harmonic maps, and show that they are solutions of nonlinear PDE system $\tau_H(f) = 0$, we also call them Euler–Lagrange equations. By the second variation formula (see Theorem \ref{thm2.4}), we prove that if the target manifold has nonpositive curvature, the exponentially subelliptic harmonic map is stable. In particular, we obtain the instability of exponentially subelliptic harmonic maps when the target manifold is a sphere.

\section{The first and second variation formula} 
Let $(M^{m+d},g)$ be a compact connected $(m + d)$ dimensional Riemannian manifold, and $g$ is the Riemannian metric on $M$. Let $H$ be a m-dimensional subbundle of tangent bundle $TM$ and set $g_H = g|H$. We call $(M, H, g_H)$ a sub-Riemannian manifold. With respect to $g$, we have the following orthogonal decomposition of $TM$:
\begin{equation}
TM = H \oplus V,
\end{equation} where $H$ (resp. $V$) is referred to as the horizontal (resp. vertical) distribution of $(M, H, g_H)$. Denote $\nabla$ the Levi-Civita connection of $(M, g)$. Then we have the induced connections $\nabla^H$ and $\nabla^V$ on $H$ and $V$ respectively \cite{CDY22}: \begin{equation}
\begin{aligned}
\nabla^H_X Y = \pi_H(\nabla_X Y), \\
\nabla^V_X Z= \pi_V(\nabla_X Z),
\end{aligned}
\end{equation}for any $X \in TM, Y \in \Gamma(H), Z\in \Gamma(V)$, where $\pi_H:TM \to H$ and $\pi_V:TM \to V$ are are the projection morphisms on $H$ and $V$ respectively. For convenience, we shall use the following convention on the range of indices on $M$:
\begin{equation} 
\begin{aligned} 
	&A,B,C,\cdots = 1, 2,\dots, m+d;\\
	&i,j,k,\dots = 1,2, \dots, m,\\
	&\alpha,\beta,\gamma,\cdots = m+1,\dots,m+d;
\end{aligned} 
\end{equation}
where we agree that repeated indices are summed over the respective ranges.

Let $\{e_A\}=\{e_i,e_\alpha\}$ be a local orthonormal frame field of $TM$, where $e_i \in \Gamma(H), e_\alpha \in \Gamma(V)$,  and let $\{\omega^A\}=\{\omega^i,\omega^\alpha\}$ be its coframe field. Let $(N^n
, h)$ be another Riemannian manifold with the Levi-Civita connection $\tilde{\nabla}$, and $f:(M, H, g_H) \to (N^n, h)$ be a smooth map from $M$ into $N$. Denote $\tilde{\nabla}^f$ the pull-back connection of $\tilde{\nabla}$ on $N$ by $f$, then we can define the following second fundamental form \cite{CDY22} with respect to $(\tilde{\nabla}, \tilde{\nabla}^f)$: \begin{equation}
\begin{aligned}
\beta_H(f)(X,Y) = \tilde{\nabla}^f_Y(df_H(X))-df_H(\nabla_Y X), \\
\beta_V(f)(X,Y) = \tilde{\nabla}^f_Y(df_V(X))-df_V(\nabla_Y X),
\end{aligned}
\end{equation}for any $X,Y \in \Gamma(TM)$, where $df_H=df \circ \pi_H$ and $df_V=df \circ \pi_V$. With respect to $\tilde{\nabla}^f$, we introduce the exterior covariant differential operator $d:\mathcal{A}^p(\xi) \to :\mathcal{A}^{p+1}(\xi) $ and its codifferential operator $d^*:\mathcal{A}^p(\xi) \to :\mathcal{A}^{p-1}(\xi)$, where $\xi = f^{-1}TN$ (cf. \cite{EL83}).

\begin{lemma}\label{lem 2.1} If $f:(M, H, g_H) \to (N^n, h)$ is a smooth map, then for any $W \in f^{-1}TN$ we have
\begin{equation} \begin{aligned}
\int_{M} e^{\frac{|df_H|^2}{2}}[e_i \langle W, df_H(e_i)\rangle-\langle W, df_H(\nabla_{e_i}e_i ) \rangle ] dv_g \\ = \int_{M} e^{\frac{|df_H|^2}{2}}\langle W, df(\zeta)-df_H(\nabla \frac{|df_H|^2}{2}) \rangle dv_g,
\end{aligned}
\end{equation}where $\zeta = \pi_H(\nabla_{e_\alpha}e_\alpha)$ and $dv_g$ is the Riemannian volume element of $g$.
\end{lemma}
\begin{proof} Define a 1-form $\theta(X) =\langle W, d_Hf(X)\rangle$ on $M$. According to the proof of \cite[Lemma 2.2]{CDY22}, we have \begin{equation}
\begin{aligned}
d^* \theta = -[e_i \langle W, df_H(e_i)\rangle-\langle W, df_H(\nabla_{e_i}e_i ) \rangle ] +\langle W, df(\zeta)\rangle.
\end{aligned}
\end{equation}
Since \begin{equation}
\int_{M} e^{\frac{|df_H|^2}{2}} div(\theta)  dv_g= \int_{M} [div(e^{\frac{|df_H|^2}{2}} \theta) - e^{\frac{|df_H|^2}{2}}\langle \theta, \nabla \frac{|df_H|^2}{2} \rangle] dv_g,
\end{equation}then we get the result by the divergence theorem.
\end{proof}

Let $f:(M, H, g_H) \to (N^n, h)$ be a smooth map from a sub-Riemannian manifold into another Riemannian manifold, we can define the horizontal exponential energy functional as \begin{equation}
E(f) = \int_{M} e^{\frac{|df_H|^2}{2}} dv_g,
\end{equation}where $e_f = \frac{1}{2}|df_H|^2$ is the horizontal energy density of $f$.

\begin{definition}
A critical map of $E(f)$ is called a exponentially subelliptic harmonic map.
\end{definition} 

\begin{theorem} \label{th2.3}
Let $\{f_t\}_{|t|<\epsilon} \in C^\infty(M \times (-\epsilon, \epsilon), N)$ be a one-parameter family of maps $f_t$ with $f_0 =f$ such that the variation vector field $V = \frac{\partial f_t}{\partial t}|_{t=0} \in \Gamma(f^{-1}TN)$ has compact support. Then we obtain the first variation formula of the horizontal exponential energy
\begin{equation}
\frac{d}{dt} E(f_t)|_{t=0} = - \int_{M} e^{\frac{|df_H|^2}{2}}\langle V, \tau_H(f) \rangle dv_g,
\end{equation} where $\tau_H(f) = \beta_H(f)(e_i,e_i) - df(\zeta)+df_H(\nabla \frac{|df_H|^2}{2})$ is called the exponential tension field of $f$. It implies that a map $f$ is exponentially subelliptic harmonic iff it satisfies the Euler-Lagrange equation $\tau_H(f)=0$.
\end{theorem} 
\begin{proof}
Define a map $F: M \times (-\epsilon, \epsilon) \to N$ by $F(x,t) = f_t(x)$. Then we obtain \begin{equation}  \begin{aligned} \label{2.10}
\frac{d}{dt} E(f_t) &= \int_{M} e^{\frac{|d{f_t}_H|^2}{2}} \langle \tilde{\nabla}_{\frac{\partial}{\partial t}} dF(e_i), dF(e_i) \rangle dv_g \\
&= \int_{M} e^{\frac{|d{f_t}_H|^2}{2}} \langle \tilde{\nabla}_{e_i} dF({\frac{\partial}{\partial t}}), dF(e_i) \rangle dv_g.
\end{aligned}
\end{equation}
Letting $t = 0$ in \eqref{2.10} and using Lemma \ref{lem 2.1}, we have
\begin{equation}
\begin{aligned}
\frac{d}{dt} E(f_t)|_{t=0} &= \int_{M} e^{\frac{|df_H|^2}{2}} \langle \tilde{\nabla}_{e_i} V, df(e_i) \rangle dv_g \\
&= \int_{M} e^{\frac{|df_H|^2}{2}} [{e_i}\langle V, df(e_i) \rangle - \langle V, \tilde{\nabla}_{e_i}df(e_i) \rangle ]dv_g \\
&= \int_{M} e^{\frac{|df_H|^2}{2}} [{e_i}\langle V, df(e_i) \rangle - \langle V, \beta_H(f)(e_i,e_i)  + df_H(\nabla_{e_i}e_i) \rangle ]dv_g  \\
&= \int_{M} e^{\frac{|df_H|^2}{2}} [- \langle V, \beta_H(f)(e_i,e_i)  \rangle+ \langle V, df(\zeta)-df_H(\nabla \frac{|df_H|^2}{2}) \rangle]dv_g 
\\
&= - \int_{M} e^{\frac{|df_H|^2}{2}} [\langle V, \beta_H(f)(e_i,e_i) - df(\zeta)+df_H(\nabla \frac{|df_H|^2}{2}) \rangle]dv_g.
\end{aligned}
\end{equation}
\end{proof}

Then the second variation of the horizontal exponential energy of $f$ is as follows
\begin{theorem} \label{thm2.4}
Let $f:(M, H, g_H) \to (N^n, h)$ be an exponentially subelliptic harmonic, then \begin{equation}
\frac{d^2}{dt^2} E(f_t)|_{t=0} = \int_{M} e^{\frac{|df_H|^2}{2}} [\langle \tilde{\nabla}_{e_i}V, df(e_i) \rangle^2 + \langle \tilde{\nabla}_{e_i}V,  \tilde{\nabla}_{e_i}V \rangle -\langle \tilde{R}(df(e_i),V)V, df(e_i) \rangle]dv_g,
\end{equation}where $\tilde{R}$ is the curvature tensor of $N$.
\end{theorem}
\begin{proof}
It follows from \eqref{2.10} that \begin{equation}
\begin{aligned} \label{2.13}
\frac{d^2}{dt^2} E(f_t) &=  \int_{M} e^{\frac{|d{f_t}_H|^2}{2}} \langle \tilde{\nabla}_{e_i} dF({\frac{\partial}{\partial t}}), dF(e_i) \rangle^2 dv_g \\ & \ \ \ \ \ + \int_{M} e^{\frac{|d{f_t}_H|^2}{2}} [\langle \tilde{\nabla}_{\frac{\partial}{\partial t}}\tilde{\nabla}_{e_i} dF({\frac{\partial}{\partial t}}), dF(e_i) \rangle + \langle \tilde{\nabla}_{e_i} dF({\frac{\partial}{\partial t}}), \tilde{\nabla}_{\frac{\partial}{\partial t}}dF(e_i) \rangle] dv_g\\
&= \int_{M} e^{\frac{|d{f_t}_H|^2}{2}} [\langle \tilde{\nabla}_{e_i} dF({\frac{\partial}{\partial t}}), dF(e_i) \rangle^2 + \langle \tilde{\nabla}_{e_i} dF({\frac{\partial}{\partial t}}), \tilde{\nabla}_{e_i}dF({\frac{\partial}{\partial t}}) \rangle]dv_g \\ 
&\ \ \ \ \ +\int_{M} e^{\frac{|d{f_t}_H|^2}{2}} [-\tilde{R}(e_i,{\frac{\partial}{\partial t}})dF({\frac{\partial}{\partial t}}), dF(e_i)+ \langle \tilde{\nabla}_{e_i}\tilde{\nabla}_{\frac{\partial}{\partial t}} dF({\frac{\partial}{\partial t}}), dF(e_i) \rangle]
dv_g\end{aligned}.
\end{equation}
Putting $t = 0$ in \eqref{2.13} and using Lemma \ref{lem 2.1} again, we have \begin{equation}
\begin{aligned}
\frac{d^2}{dt^2} E(f_t)|_{t=0} &= \int_{M} e^{\frac{|df_H|^2}{2}} [\langle \tilde{\nabla}_{e_i}V, df(e_i) \rangle^2 + \langle \tilde{\nabla}_{e_i}V,  \tilde{\nabla}_{e_i}V \rangle-\langle \tilde{R}(df(e_i),V)V, df(e_i) \rangle \\
&\ \ \ \ \ \ \ \ \ \ +{e_i}\langle \tilde{\nabla}_{\frac{\partial}{\partial t}}V, df(e_i) \rangle - \langle \tilde{\nabla}_{\frac{\partial}{\partial t}}V, \tilde{\nabla}_{e_i}df(e_i) \rangle ]dv_g \\
&= \int_{M} e^{\frac{|df_H|^2}{2}} [\langle \tilde{\nabla}_{e_i}V, df(e_i) \rangle^2 + \langle \tilde{\nabla}_{e_i}V,  \tilde{\nabla}_{e_i}V \rangle-\langle \tilde{R}(df(e_i),V)V, df(e_i) \rangle \\
&\ \ \ \ \ \ \ \ \ \ +{e_i}\langle \tilde{\nabla}_{\frac{\partial}{\partial t}}V, df(e_i) \rangle - \langle \tilde{\nabla}_{\frac{\partial}{\partial t}}V, \beta_H(f)(e_i,e_i)  + df_H(\nabla_{e_i}e_i) \rangle ]dv_g \\
&= \int_{M} e^{\frac{|df_H|^2}{2}} [\langle \tilde{\nabla}_{e_i}V, df(e_i) \rangle^2 + \langle \tilde{\nabla}_{e_i}V,  \tilde{\nabla}_{e_i}V \rangle-\langle \tilde{R}(df(e_i),V)V, df(e_i) \rangle \\
&\ \ \ \ \ \ \ \ \ \ - \langle \tilde{\nabla}_{\frac{\partial}{\partial t}}V, \beta_H(f)(e_i,e_i) - df(\zeta)+df_H(\nabla \frac{|df_H|^2}{2}) \rangle]dv_g \\
&=\int_{M} e^{\frac{|df_H|^2}{2}}[\langle \tilde{\nabla}_{e_i}V, df(e_i) \rangle^2 + \langle \tilde{\nabla}_{e_i}V,  \tilde{\nabla}_{e_i}V \rangle -\langle \tilde{R}(df(e_i),V)V, df(e_i) \rangle]dv_g,
\end{aligned}
\end{equation}where the last equality is obtained by the Euler-Lagrange equation.
\end{proof}

\section{Stability of exponentially subelliptic harmonic maps}
Suppose that $f:(M, H, g_H) \to (N^n, h)$ is an exponentially subelliptic harmonic map. We put \begin{equation}
I_H(V,V) = \int_{M} e^{\frac{|df_H|^2}{2}}[\langle \tilde{\nabla}_{e_i}V, df(e_i) \rangle^2 + \langle \tilde{\nabla}_{e_i}V,  \tilde{\nabla}_{e_i}V \rangle -\langle \tilde{R}(df(e_i),V)V, df(e_i) \rangle]dv_g.
\end{equation}

\begin{definition}
An exponentially subelliptic harmonic map is called stable if $I_H(V,V) \ge 0$ for any compactly supported vector field $V$ along $f$.
\end{definition}

\begin{proposition}
If the curvature of target manifold $N$ is always nonpositive, then $f$ is stable.
\end{proposition}

Next, we consider the stability problem when the target manifold is a sphere $S^n$ as a submanifold in $\mathbb{R}^{n+1}$. Let $^{S}\tilde{\nabla}$ and $^{R}\tilde{\nabla}$ be the Levi-Civita connections on $S^n$ and $\mathbb{R}^{n+1}$, respectively. For the vector $V$ in $\mathbb{R}^{n+1}$ at $x \in S^n$, we have the following orthogonal decomposition:\begin{equation}
V = V^T + V^{\perp},
\end{equation} where $V^T$ is tangent to $S^n$ and
$V^{\perp}=\langle V, x\rangle x$ is the normal part to $S^n$. Let $B$ denote the second fundamental form of $S^n$ in $\mathbb{R}^{n+1}$. For tangent vector fields $X$, $Y$ and normal vector field $W$ of $S^n$ around $x$, the self-adjoint operator $A^W$ of $W$ is defined by \begin{equation}
A^W(X) = -(^{R}\tilde{\nabla}_X W)^T, 
\end{equation}it is known that \begin{equation}
\langle A^W(X), Y\rangle = \langle B(X,Y), W\rangle =- \langle X, Y \rangle \langle x,W \rangle.
\end{equation}

\begin{theorem}
Suppose $f:M \to S^n$ is a nonconstant exponentially subelliptic harmonic map from a compact sub-Riemannian manifold $M$ into the unit sphere $S^n$. If 
\begin{equation}
\int_{M} e^{\frac{|df_H|^2}{2}}[|df_H|^2(|df_H|^2-(n-2))]dv_g <0,
\end{equation}then $f$ is unstable.
\end{theorem}

\begin{proof}
Let $\tilde{\nabla}$ and ${\tilde{R}^S}$ denote the induced connection on $f^{-1}TS^n$ and the curvature tensor of $S^n$, respectively. Choose a parallel orthonormal frame field $\{V_a\}_{a=1}^{n+1}$ in $\mathbb{R}^{n+1}$. Then from Theorem \ref{thm2.4} that the second variation is
\begin{equation}\begin{aligned}
I_H(V_a^T,V_a^T) &= \int_{M} e^{\frac{|df_H|^2}{2}}[\langle \tilde{\nabla}_{e_i}V_a^T, df(e_i) \rangle^2 + \langle \tilde{\nabla}_{e_i}V_a^T,  \tilde{\nabla}_{e_i}V_a^T \rangle \\ & \ \ \ \ \ \ \ \  -\langle \tilde{R}^S(df(e_i),V_a^T)V_a^T, df(e_i) \rangle]dv_g.
\end{aligned}
\end{equation}From $\langle \tilde{\nabla}_{e_i}V_a^T, df(e_i) \rangle = -|df(e_i)|^2 \langle x,V_a \rangle $ and $\langle \tilde{\nabla}_{e_i}V_a^T,  \tilde{\nabla}_{e_i}V_a^T \rangle = |df(e_i)|^2 \langle x,V_a \rangle^2$ (cf. \cite{Ch16}), it follows that \begin{equation} \label{3.7}
\sum_{a}[\sum_{i}\langle \tilde{\nabla}_{e_i}V_a^T, df(e_i) \rangle]^2 = |df_H|^4, \end{equation}
\begin{equation} \label{3.8}
\sum_{a}\sum_{i}\langle \tilde{\nabla}_{e_i}V_a^T,  \tilde{\nabla}_{e_i}V_a^T \rangle  = |df_H|^2.
\end{equation}Furthermore, \begin{equation}\begin{aligned}
\sum_{a}\sum_{i}\langle \tilde{R}^S(df(e_i),V_a^T)V_a^T, df(e_i) \rangle &= \sum_{a}\sum_{i}|df(e_i)|^2|V_a^T|^2-\langle df(e_i), V_a^T \rangle^2 \\ &= |df_H|^2 \sum_{a}|V_a^T|^2 - |df_H|^2\\ &=(n-1)|df_H|^2. \label{3.9}
\end{aligned}\end{equation}
Combining \eqref{3.7}, \eqref{3.8} and \eqref{3.9} yields \begin{equation} \begin{aligned} \label{3.10}
\sum_{a}I_H(V_a^T,V_a^T) &= \int_{M} e^{\frac{|df_H|^2}{2}}[|df_H|^4+ |df_H|^2-(n-1)|df_H|^2]dv_g \\
&= \int_{M} e^{\frac{|df_H|^2}{2}}[|df_H|^2(|df_H|^2+(2-n))]dv_g .
\end{aligned} \end{equation}
By \eqref{3.10} and the assumption conditions, we have \begin{equation}
\sum_{a}I_H(V_a^T,V_a^T)<0,
\end{equation}which implies $f$ is unstable.
\end{proof}
\begin{corollary}
If $|df_H|^2<n-2$, then any nonconstant exponentially subelliptic harmonic map from a compact sub-Riemannian manifold $M$ into $S^n$ is unstable.
\end{corollary}

\bigskip
Xin Huang

School of Mathematical Sciences

Nanjing University of Information Science and Technology

Nanjing 210044, P. R. China

003941@nuist.edu.cn

    \end{document}